\documentclass{amsart}
\usepackage[all]{xy}
\usepackage{amsmath,amssymb}

\theoremstyle{plain}

\newtheorem{theorem}{Theorem}[section]

\newtheorem{corollary}[theorem]{Corollary}

\theoremstyle{definition}

\newtheorem{example}[theorem]{Example}

\theoremstyle{remark}
\newtheorem{remark}[theorem]{Remark}

\newcommand{\secref}[1]{Section~\ref{#1}}
\newcommand{\thmref}[1]{Theorem~\ref{#1}}

\newcommand{\exref}[1]{Example~\ref{#1}}

\newcommand{\remref}[1]{Remark~\ref{#1}}

\def\implies{\Longrightarrow}

\def\lto{\longrightarrow}

\def\Z{{\mathbb Z}}
\def\Q{{\mathbb Q}}
\def\C{{\mathbb C}}

\def\map{\mathrm{map}}
\def\Baut{\mathrm{B \, aut}}
\def\Der{\mathrm{Der}}
\def\aut{\mathrm{aut}}
\def\Hom{\mathrm{Hom}}

\def\ker{\mathrm{ker}}
\def\im{\mathrm{im}}
\def\cat{\mathrm{cat}}
\def\dim{\mathrm{dim}}

\begin{document}

\title[The Evaluation Subgroup of a Fibre Inclusion]
{The Evaluation Subgroup of a Fibre Inclusion}

\author{Gregory  Lupton}

\address{Department of Mathematics,
          Cleveland State University,
          Cleveland OH 44115}

\email{G.Lupton@csuohio.edu}

\author{Samuel Bruce Smith}

\address{Department of Mathematics,
  Saint Joseph's University,
  Philadelphia, PA 19131}

\email{smith@sju.edu}

\date{\today}



 \begin{abstract}
 Let  $\xi\colon X \stackrel{j}{\rightarrow} E \stackrel{p}{\rightarrow}
 B$ be a fibration of  simply connected CW complexes of finite type
 with classifying map $h \colon B \to \Baut_1(X).$
We study  the evaluation subgroup $G_n(E, X;j)$ of the fibre
inclusion as an invariant of the fibre-homotopy type of $\xi$.
For spherical fibrations, we show the evaluation subgroup may be
expressed as
 an extension of the  Gottlieb group of
the fibre sphere  provided the classifying map $h$ induces the trivial
map on homotopy groups. 
We extend this result after rationalization: We show that  
 the decomposition $G_*(E, X; j) \otimes \Q = (G_*(X) 
\otimes \Q) \oplus (\pi_*(B) \otimes \Q)$  is equivalent to
the condition $(h_\sharp)_\Q = 0$. 
   \end{abstract}
\maketitle

\section{Introduction.}
The $n$th Gottlieb group $G_n(X)$ of a space $X$ is the subgroup of
$\pi_n(X)$ consisting of homotopy classes of maps $g: S^n \to X$
such that the map $(g\mid 1_X) \colon S^n \vee X \to X$, defined on
the wedge, extends to some map $F \colon S^n \times X \to X$  of the
product. Alternately, $G_n(X)$ is the image of the map induced on 
homotopy groups by the evaluation map $\omega \colon \aut_1(X) 
\to X,$ where $\aut_1(X)$  is the space of continuous functions 
homotopic to the identity    of $X$.  
 The definition may be generalized by replacing the
identity  by an arbitrary based map $f \colon X
\to Y$: the \emph{$n$th evaluation subgroup $G_n(Y, X;f)$ of the map
$f$} is the subgroup of $\pi_n(Y)$ represented by
maps $g \colon S^n \to Y$ such that  $(g\mid f) \colon S^n \vee X
\to Y$ extends to some map $F\colon S^n \times X \to Y$. The
generalization provides some functorality in that  the map on homotopy
groups induced by $f$ restricts to a map $f_\sharp
\colon G_n(X) \to G_n(Y,X ;f )$.  On the other hand, it is
well-known that the Gottlieb group fails to be a functor since,
generally, a map $f \colon X \to Y$ does not yield a
homomorphism $f_\# \colon G_n(X) \to G_n(Y)$.

The Gottlieb groups $G_n(X)$ play a well-known role in the homotopy
theory of fibrations with fibre a  CW complex $X$ of finite type.
Notably,  the  result of Gottlieb \cite[Th.2]{Got2} identifies
$G_*(X)$ as the image of the connecting homomorphism in the long exact
homotopy sequence of the universal fibration with fibre $X$ (see 
diagram (\ref{eq:connect}) below).  In this
paper, we investigate the role of the more general evaluation
subgroups for the homotopy theory of fibrations.

Let $\xi\colon X \stackrel{j}{\rightarrow} E
\stackrel{p}{\rightarrow} B$ be a fibration of  simply connected CW
complexes of finite type. We  consider  the homomorphism $j_\sharp
\colon G_n(X) \to G_n(E, X;j)$ induced by the fibre inclusion of
$\xi$.  By the universal property of the Gottlieb group mentioned
above, restricting the homomorphisms in the long exact homotopy
sequence of $\xi$ yields a   sequence of groups and homomorphisms
\begin{equation}
\label{eq:Gottlieb sequence}
 \xymatrix{ \cdots \ar[r] & \pi_{n+1}(B) \ar[r]^{\partial}
 &   G_n(X) \ar[r]^{\! \!  \! \! \! \! \! \! j_\sharp} & G_n(E, X;j) \ar[r]^{\ \ p_\sharp}
 &  \pi_n(B) \ar[r] & \cdots}.
\end{equation}
We call (\ref{eq:Gottlieb sequence}) the \emph{Gottlieb sequence of
the fibration $\xi$}.  In general, it is not exact. 
 Our overriding purpose in this paper is to show that
exactness properties of the Gottlieb sequence represent an
interesting measure of the relative triviality of the fibration.

Successive compositions of maps in the Gottlieb sequence are zero,
since they are restrictions of those in the long exact homotopy
sequence of $\xi$. Thus we may consider homology groups at each type
of term. We focus on the generalized evaluation  subgroup terms and define the
$n$th \emph{Gottlieb homology group of $\xi$} to be the subquotient
$$GH_n(\xi) = \frac{\ker\{ p_\sharp \colon G_n(E, X;j) \to \pi_n(B) \} }
{\im \{ j_\sharp \colon G_n(X) \to G_n(E, X;j) \}}
$$
of  $G_n(E, X;j)\subseteq \pi_n(E)$.  As we
shall see, the Gottlieb homology detects the nontriviality of
the fibration in the sense of fibre-homotopy type. We  say the
fibration $\xi$ is {\em Gottlieb trivial in degree $n$} if (\ref{eq:Gottlieb
sequence}) breaks into a short exact sequence
\begin{equation}
\label{eq:G-seq}
\xymatrix{ 0 \ar[r] & G_n(X) \ar[r]^-{j_\sharp} & G_n(E, X;j)
\ar[r]^-{p_\sharp} & \pi_n(B) \ar[r] & 0}
\end{equation}
in degree $n,$  and {\em Gottlieb trivial} if this occurs in all 
degrees $n \geq 2$.  

Recall that a fibration $\xi$ as above is classified by
a map  $h \colon B \to \Baut_1(X),$ where   
 $\Baut_1(X)$ is the classifying
space of the monoid  $\aut_1(X)$ in the sense of Dold-Lashof \cite{D-L} (see \cite{All, St,
D}).   Our main results relate
the notion of Gottlieb triviality to the vanishing of the
homomorphism induced on
homotopy groups  by the classifying map of the fibration. 
In \secref{sec:spherical}, we make
 some observations in ordinary homotopy
theory that relate the behaviour of the Gottlieb sequence with that
of the classifying map.   We then prove that a
spherical fibration is Gottlieb trivial whenever $h_\sharp = 0$
(\thmref{thm:GH=0}).     

In \secref{sec:rational}, we
obtain a complete result in this vein within the framework of rational homotopy theory. 
We prove  that, after tensoring with the rationals,   the Gottlieb sequence of a  
fibration with fibre an arbitrary finite complex
splits as in (\ref{eq:G-seq})   if and only if the homomorphism induced on rational homotopy groups 
by the classifying map  is trivial 
(\thmref{thm:rational split}).  
This result depends  upon a
 description of the map $(h_\sharp)_\Q$  in terms of 
derivations of the Sullivan model of the fibre. We give this 
description in \secref{sec:h}.   The 
key step here comes from a close study of the holonomy action of a 
fibration. This  description of $(h_\sharp)_\Q$  (\thmref{thm:h})
 fits the classifying map into the framework developed
in \cite{L-S1}, whereby chain complexes of (generalized) derivations
of Sullivan models are used to describe the rationalized evaluation
subgroups. 

\begin{remark}\label{rem:Lee-Woo}
The Gottlieb sequence of a fibration $\xi : F \stackrel{j}{\to} E
\stackrel{p}{\to} B$ is a special case of the $G$-sequence of a
map  as developed and studied by Lee and Woo
(cf. \cite{L-W4}). 
We have elected to focus on the Gottlieb sequence of a fibration
since the construction is direct  and avoids
consideration of the  relative term which occurs in the general
$G$-sequence.   However, the results of this paper may generally  be
rephrased in terms of the $G$-sequence and the $\omega$-homology of
a map.  In particular, our definition of the Gottlieb homology
corresponds to that of the $\omega$-homology group $H_*^{\omega
b}(E, X; j)$ \cite{L-W4}.
\end{remark}

\section{First Results in Ordinary Homotopy Theory} 
\label{sec:spherical}
Here and throughout, all spaces  will be assumed to be based
 CW complexes of finite type.   A fibration  $\xi
\colon X \stackrel{j}{\to} E \stackrel{p}{\to} B$ will be  a
Hurewicz fibration \cite[p.66]{Sp}. The map $j$ will usually denote
the inclusion of the fibre $X  = p^{-1}(b_0).$  However, by
\cite[Prop.1.4]{Go1}, the
evaluation subgroups are homotopy invariants of a map and so we may, when
necessary, assume $j$ is only a fixed equivalence from $X$ to the
actual fibre over the basepoint. We will write $\partial \colon \Omega 
B \to X$ for the connecting map in the Puppe sequence for $\xi$ 
\cite[Th.III.6.22]{Wh}.  After the identification $\pi_n(\Omega B) 
\cong \pi_{n+1}(B),$ $\partial$ induces    the connecting homorphism
$\partial_\sharp \colon \pi_{n+1}(B) \to \pi_n(X)$ in the long exact 
homotopy sequence for $\xi$.  
We let  $\partial_\infty \colon \Omega \Baut_1(X) \to X$ and 
$\partial_U \colon \Omega \Baut_1(X) \to   \aut_1(X)$ denote the 
connecting maps in the universal $X$-fibration
$X \to UX \to \Baut_1(X)$    and the universal principal 
$\aut_1(X)$-fibration $\aut_1(X) \to U\aut_1(X) \to \Baut_1(X),$ 
respectively. 
By results of Gottlieb \cite[\S 4]{Go1} and Dror-Zabrodsky 
 \cite[Prop.4.1]{D-Z}, we then have the following  homotopy commutative diagram 
\begin{equation} \label{eq:connect}
\xymatrix{ \Omega B \ar[drr]_{\partial} \ar[rr]^{\Omega h} && \Omega 
\Baut_1(X) \ar[d]^{\partial_\infty} \ar[rr]_{\simeq}^{\partial_U} && \aut_1(X) 
\ar[dll]^{\omega} \\
& & X. & }
\end{equation}

The  Gottlieb sequence is an invariant of the  fibre-homotopy 
type  of $\xi$. Precisely, say two group homomorphisms $\phi \colon
G_1 \to G_2$ and $\psi \colon H_1 \to H_2$  are {\em equivalent} if
there are isomorphisms $\theta_i \colon G_i \to H_i$   such that
$\psi \circ \theta_1 = \theta_2 \circ \phi.$ This notion of
equivalence extends in an obvious way to any sequence or commutative
diagram of groups and homomorphisms. Suppose $\xi_1 \colon X_1
\stackrel{j_1}{\to} E_1 \stackrel{p_1}{\to} B$ is fibre-homotopy
equivalent to  $\xi_2 \colon X_2 \stackrel{j_2}{\to} E_2
\stackrel{p_2}{\to} B$ \cite[p.100]{Sp}. Using \cite[Prop.1.4]{Go1} 
again, we obtain that
the equivalence between long exact homotopy sequences induced by the
fibre-homotopy equivalence preserves the evaluation subgroups and
thus restricts to an  equivalence of Gottlieb sequences.

We next observe that Gottlieb triviality fits properly between two
familiar notions of the relative triviality of a fibration. Say 
$\xi \colon X \stackrel{j}{\to} E \stackrel{p}{\to} B$ is {\em
fibre-homotopically trivial} if $\xi$ is fibre-homotopy equivalent  to
the product
 $\pi \colon X
\stackrel{i_2}{\to} B \times X \stackrel{p_1}{\to} B.$ Say $\xi$
is {\em weak-homotopically  trivial} if 
$\partial_\# = 0   \colon \pi_{n+1}(B) \to \pi_n(X)$   for all $n.$ 
\begin{theorem}\label{thm:implications}
For any fibration we have the following implications:
$$ \hbox{fibre-homotopically trivial }  \implies \hbox{Gottlieb
trivial  } \implies \hbox{weak-homotopically trivial}.$$
Furthermore, each of the reverse implications fails.
\end{theorem}

\begin{proof}
An easy direct argument shows that a trivial fibration is Gottlieb 
trivial (cf.  \cite[Cor.13]{L-W1}).   The first implication follows 
from this and the 
fibre-homotopy invariance mentioned above. For the second implication, note that $p_\sharp$
surjective implies $\partial_\# = 0$. A separating example for the 
first implication is given by \exref{ex:principal}. Separating 
examples for the second are given by 
Examples~\ref{ex:integral not Gottlieb trivial} and \ref{ex:not Gottlieb}.
\end{proof}

\begin{example}
\label{ex:principal}
Let  $\xi: G \stackrel{j}{\to} E \stackrel{p}{\to} B$
be  a principal bundle with structure group $G$.
It is well-known
that the Gottlieb groups of $G$ coincide with the homotopy groups of
$G$.  This identity extends to the evaluation subgroups $G_n(E, G;j)$:
For   let $\alpha : S^n \to E$ represent an arbitrary
homotopy class. Following Steenrod \cite[\S 8.7]{Steen},   consider the   principal map
  $P : G \times E \to E$   induced by the action of $G$ on the fibres of
  $p : E \to B.$
Define $F \colon S^n \times G \to
E$  by $F(s, g) = P(g, \alpha(s)).$ It is easy to check $F$ extends
$(\alpha \, | \, j) \colon S^n \vee G \to E$ to the product and the
claim follows.
Thus for principal bundles  the Gottlieb sequence
coincides with the long exact homotopy sequence and we have   Gottlieb trivial
is equivalent to  weak-homotopically trivial.
The  above analysis holds as well for principal $H$-fibrations
    as in \cite{D-L}.
We thus obtain Gottlieb trivial fibrations which are not
fibre-homotopy trivial by considering, for example,  the (nontrivial)
stages  in the Postnikov
decomposition of a space.
\end{example}

Given a fibration $\xi \colon X \stackrel{j}{\to} E \stackrel{p}{\to} 
B$ with classifying map
$h\colon B
\to \Baut_1(X)$ recall that $\xi$ is 
 fibre-homotopically trivial if and only if $h$ is null-homotopic. 
 The weaker condition 
$h_\sharp = 0$    is thus an approximation to $\xi$ 
being fibre-homotopically trivial. In what follows we compare this 
condition to Gottlieb triviality. We begin with the following general fact.

\begin{theorem}\label{thm:p onto ker}
Let $\xi \colon X \stackrel{j}{\to} E \stackrel{p}{\to} B$ be a
fibration of simply connected, finite type CW complexes with
classifying map $h\colon B \to \Baut_1(X)$. Then we have
$$p_\sharp(G_n(E, X;j)) \supseteq \ker\{h_\sharp \colon \pi_n(B) \to
\pi_n(\Baut_1(X)) \}.$$
\end{theorem}

\begin{proof} Let $\alpha \colon S^n \to B$ represent a homotopy class
in $\pi_n(B)$ with $h_\sharp(\alpha) = 0.$ Let $X
\stackrel{j^*}{\to} E^* \stackrel{p^*}{\to} S^n$ denote the
pullback of $X \stackrel{j}{\to} E \stackrel{p}{\to} B$ by $\alpha
\colon S^n \to B$ and $\alpha^* \colon E^* \to E$ the fibre map
covering $\alpha.$ By hypothesis, the composition $h \circ  \alpha
: S^n \to \Baut_1(X)$ is homotopically trivial. This means $X
\stackrel{j^*}{\to} E^* \stackrel{p^*}{\to} S^n$ is 
fibre-homotopically trivial. We thus have a fibre-homotopy equivalence
$H \colon S^n \times X \to E^*$ giving a commutative diagram
$$\xymatrix{ &  X \ar[dl]_{i_2} \ar[d]_{j^*} \ar[dr]^{j} \\
S^n\times X \ar[r]^{\ \ H} \ar[dr]_{\pi_1}  & E^*
\ar[r]^{\alpha^*}
\ar[d]_{\! \! \! \! p^*} & E \ar[d]^{p} \\
& S^n \ar@/^2pc/[ul]_{i_1} \ar[r]^{\alpha}
\ar@{.>}[ur]^{\tilde{\alpha}} & B.}
$$
Define a lifting $\tilde{\alpha} \colon S^n \to E$ by
$\tilde{\alpha} = \alpha^* \circ H \circ i_1.$ To see
$\tilde{\alpha}$ represents a class in $G_n(E, X; j)$ define    $G
\colon S^n \times X  \to E$ by   $G = \alpha^* \circ H$ and observe that
$G$ extends $(\tilde{\alpha}  \mid  j ) \colon S^n \vee X \to
E$ to the product.
\end{proof}

As a consequence, we obtain: 

\begin{theorem}\label{thm:one-to-one}
Let $\xi \colon X \stackrel{j}{\to} E \stackrel{p}{\to} B$ be a
fibration of simply connected, finite type CW complexes with
classifying map $h\colon B \to \Baut_1(X)$. 
Suppose $ h_\sharp = 0 \colon \pi_n(B) \to \pi_n(\Baut_1(X))$
for some $n \geq 2.$   Then 
$j_\# \colon G_{n-1}(X) \to
G_{n-1}(E, X; j)$  is injective and 
  $p_\# \colon G_n(E, X; j) \to \pi_n(B)$  is surjective. 
\end{theorem}
\begin{proof}
Our hypothesis and  the factorization $\partial = \partial_\infty \circ \Omega h
$ given by  (\ref{eq:connect}) imply   $j_\sharp \colon 
\pi_{n-1}(X) \to \pi_{n-1}(E)$ is injective and so the restriction of 
$j_\sharp$ to $G_n(X)$ is also. 
Surjectivity  is direct from  \thmref{thm:p onto ker}.
\end{proof}

Focusing on the case of a spherical fibration, we obtain  
the following result:
\begin{theorem}\label{thm:GH=0}
Let $\xi \colon S^k \stackrel{j}{\to} E \stackrel{p}{\to} B$ be a 
fibration
with classifying map $h\colon B \to
\Baut_1(S^k)$.
Suppose $h_\sharp = 0\colon \pi_{i}(B) \to
\pi_{i}\big(\Baut_1(S^k)\big)$ for $i = n, n+1$ and $n+k.$  Then
the fibration $\xi$ is Gottlieb trivial in degree $n$,  that is, 
we have a short exact sequence
$$\xymatrix{ 0 \ar[r] & G_n(S^k) \ar[r]^-{j_\sharp} & G_n(E, S^k;j)
\ar[r]^-{p_\sharp} & \pi_n(B) \ar[r] & 0}.$$
\end{theorem}

\begin{proof}
By \thmref{thm:one-to-one},  it suffices to show  $GH_n(\xi) = 0.$ 
Let $\beta \colon S^n \to E$ represent a homotopy class in
$\ker\{ p_\sharp \colon G_n(E, S^k;j) \to \pi_n(B) \}$. Because
$\beta$ represents an element of $G_n(E, S^k;j)$, we have a map $G
\colon S^n \times S^k \to E$ that extends $(\beta \mid j) \colon
S^n \vee S^k \to E$.   Since $p\circ\beta \sim * \colon S^n \to
B$, and since $\xi$ is a fibration, we may choose a class $\alpha
\colon S^n \to S^k$ such that $j \circ \alpha \sim \beta$. We need
to show that $\alpha \in G_n(S^k)$, that is, we must produce a map
$F \colon S^n \times S^k \to S^k$ that extends $(\alpha \mid 1)
\colon S^n \vee S^k \to S^k$.  Let $\eta \colon S^{n+k-1} \to S^n
\vee S^k$ denote the Whitehead product map $\eta = [\iota_n,
\iota_k]$, which gives the cofibre sequence $S^{n+k-1}
\stackrel{\eta}{\to} S^n \vee S^k \to S^n \times S^k$.  We have a
commutative diagram
$$\xymatrix{ S^{n+k-1} \ar[d]_{\eta}   \\
S^n \vee S^k \ar[r]^-{(\alpha\mid 1)} \ar[d]_{J}  & S^k \ar[d]^{j} \\
S^n \times S^k \ar [r]_-{G} \ar@{.>}[ur]^{F} & E}$$
in which $j \circ (\alpha\mid 1)\circ \eta \sim G\circ J\circ \eta 
\sim
*$.  As observed in the proof of \thmref{thm:one-to-one}, 
we have $j_\# \colon \pi_{n+k-1}(S^k) \to \pi_{n+k-1}(E)$ injective. It
follows that $(\alpha\mid 1)\circ \eta \sim *\colon S^{n+k-1} \to
S^k$ and hence that there exists an extension of $(\alpha\mid 1)$
to a map $F\colon S^n\times S^k \to S^k$ as desired.
\end{proof}

The following example is 
primarily intended to separate Gottlieb trivial from weak-homotopically 
trivial.  It also illustrates that, at least under rather restricted 
circumstances, converses of Theorems  \ref{thm:p onto ker}, 
\ref{thm:one-to-one} and
\ref{thm:GH=0} may hold.

\begin{example} \label{ex:integral not Gottlieb trivial}
We claim that there is a fibration $\xi \colon S^3
\stackrel{j}{\rightarrow} E \stackrel{p}{\rightarrow} S^3$ that is
weak-homotopically trival but not Gottlieb trivial.  Specifically,
$p_\# \colon G_3(E, S^3;j) \to \pi_3(S^3)$ will not be surjective in
our example. First observe, that any fibration $\xi \colon S^n
\stackrel{j}{\rightarrow} E \stackrel{p}{\rightarrow} S^n$ must be
weak-homotopically trivial (at least).  Indeed, it must have a
section $\sigma \colon S^n \to E$.  This follows, since the
connecting homomorphism $\partial_\# \colon \pi_n(S^n) \to
\pi_{n-1}(S^n) = 0$ is necessarily trivial, and so $p_\# \colon
\pi_n(E) \to \pi_n(S^n)$ is surjective.  Choose $\sigma$ to be any
element of $\pi_n(E)$ for which $p_\#(\sigma) = 1 \in \pi_n(S^n)$.

Next, any fibration $\xi \colon S^n \stackrel{j}{\rightarrow} E
\stackrel{p}{\rightarrow} S^n$ satisfies
$$p_\sharp(G_n(E, S^n;j)) \subseteq \ker\{h_\sharp \colon \pi_n(S^n) \to
\pi_n(\Baut_1(S^n)) \}.$$
Together with \thmref{thm:p onto ker}, this shows that $p_\sharp(G_n(E,
S^n;j)) = \ker\{h_\sharp \}$ in this case.  For suppose that we have
$\alpha = p_\#(\beta)$ for $\alpha \in \pi_n(S^n)$ and $\beta \in
G_n(E, S^n;j)$.  Then there exists some $G \colon S^n \times S^n \to
S^n$ that extends $(\beta\mid j) \colon S^n \vee S^n \to S^n$.  We now
argue that, without loss of generality, we may assume that $p\circ G
= \alpha\circ p_1 \colon S^n \times S^n \to S^n$, where $p_1 \colon
S^n \times S^n \to S^n$ denotes projection onto the first factor.
The cofibration sequence $S^{2n-1} \to S^n \vee S^n \to S^n \times
S^n$ gives rise to a diagram of Puppe sequences
$$\xymatrix{ [S^{2n}, E] \ar[d]^{p_*} \ar[r] & [S^n \times S^n, E] \ar[d]^{p_*} \ar[r]
& [S^n \vee S^n, E] \ar[d]^{p_*} \\
[S^{2n}, S^n]\ar@/^1pc/[u]^{\sigma_*} \ar[r] & [S^n \times S^n, S^n]
\ar[r] & [S^n \vee S^n, S^n], }$$
in which the left-hand terms act on the middle terms in the usual
way.  Since $p\circ G$ and $\alpha\circ p_1$ both map to the same
element in the lower right set, we have $\alpha\circ p_1 = (p\circ
G)^{\gamma}$ for some $\gamma \in \pi_{2n}(S^n)$.  Since we have the
section $\sigma$, we may write $\gamma = p_* \sigma_*(\gamma)$ and
so we have $\alpha\circ p_1 = (p\circ G)^{p_* \sigma_*(\gamma)} =
p_*\big( G^{\sigma_*(\gamma)}\big)$.  That is, we may replace $G$ by
$G^{\sigma_*(\gamma)}$ to obtain a map that extends $(\beta\mid j)
\colon S^n \vee S^n \to S^n$ and also projects under $p$ to
$\alpha\circ p_1$.  Finally, consider the pullback of the fibration
$\xi$ over $\alpha$.  This gives a fibration $\xi^* \colon S^n
\stackrel{j^*}{\rightarrow} E^* \stackrel{p^*}{\rightarrow} S^n$
with classifying map $h\circ\alpha\colon S^n \to \Baut_1(S^n)$.  The
maps $G \colon S^n \times S^n \to E$ and $p_1 \colon S^n \times S^n
\to S^n$ that satisfy $p\circ G = \alpha\circ p_1$ define a map
$f\colon S^n \times S^n \to E^*$ into the pullback.  This map gives
a commutative diagram
$$\xymatrix{S^n \ar[d]_{i_2} \ar[r]^{1} & S^n \ar[d]^{j^*} \\
S^n \times S^n \ar[d]_{p_1} \ar[r]^{f} & E^* \ar[d]^{p^*}\\
S^n \ar[r]^{1} & S^n. }$$
Since this displays the induced fibration $\xi^*$ as
fibre-homotopically trivial, it follows that its classifying map
$h\circ\alpha$ is null-homotopic.

To complete our example, it remains to identify a specific instance
in which $h_\# \colon \pi_n(S^n) \to \pi_n\big(\Baut_1(S^n)\big)$
may be chosen non-zero.  For this, take $n = 3$. We have
$\pi_3\big(\Baut_1(S^3)\big) \cong \pi_2\big(\Omega\Baut_1(S^3)\big)
\cong \pi_2\big(\map(S^3, S^3;1)\big)$.  Since $S^3$ is an
$H$-space, the evaluation fibration $\map_*(S^3, S^3;1) \to
\map(S^3, S^3;1) \to S^3$ admits a section and it follows that
$\pi_2\big(\map(S^3, S^3;1)\big) \cong \pi_2\big(\map_*(S^3,
S^3;1)\big)$.  Using again that $S^3$ is an $H$-space, we have a
well-known homotopy equivalence of components $\map_*(S^3, S^3;1)
\simeq \map_*(S^3, S^3;0)$ and it now follows by standard methods
that $\pi_2\big(\map_*(S^3, S^3;0)\big) \cong \pi_5(S^3) \cong
\Z_2$.  In summary, we have computed that
$\pi_3\big(\Baut_1(S^3)\big) \cong \Z_2$.  Choose $h \colon S^3 \to
\Baut_1(S^3)$ to represent the non-trivial element.  This is the
classifying map of a weak-homotopically trivial fibration $S^3 \to E
\to S^3$ that is not Gottlieb trivial, as claimed.
\end{example}

By \cite[Cor.2.2]{LMW},    $G_n(S^2) \cong  \pi_n(S^3)$ for all $n$.
We complement this result with the following:    

\begin{example}
Let $\eta_2 \colon S^3  \to S^2$ denote the Hopf map. 
We claim 
$G_n(S^2, S^3; \eta_2) = \pi_n(S^2)$ for all $n$.
Write $k \colon S^2 \to 
BS^1 = K(\Z, 2)$ for  the classifying map for $\eta_2$ viewed as a  principal 
$S^1$-fibration.  Converting  $k$ to a fibration, we obtain an 
$S^3$-fibration
  $\xi: S^3 \stackrel{\eta_2}{\to} S^2 
\stackrel{k}{\to} K(\Z, 2)$.   By \thmref{thm:GH=0},  $\xi$ is Gottlieb 
trivial  in degrees $n > 2$.  We  check directly that $\xi$ is 
Gottlieb trivial in degree $2$ as well. For  note that $G_2(S^3) = 0,$
while the Whitehead identity $[\iota_2, \eta_2] = 0$ implies $ G_2(S^2, 
S^3; 
\eta_2)  = \pi_2(S^2)$ and so $k_\sharp \colon G_2(S^2, S^3; j) \to
\pi_2(K(\Z, 2))$ is an isomorphism.  
\end{example}

\section{Derivations of Sullivan Models, the Holonomy Action and  the  Classifying Map}%
\label{sec:h} 
In this section,  we   describe the map induced  on 
rational homotopy groups by a 
classifying map for a fibration, in terms of certain chain complexes of derivations of 
Sullivan models. 
We first introduce some notation   for working
in Sullivan's differential graded (DG) algebra framework for
rational homotopy theory, 
for which our general reference is \cite{F-H-T}.

By a {\em DG algebra} we  mean a pair
$A, d$ where $A$ is a  connected, commutative   graded algebra over $\Q$.  The differential  $d$
increases degree by one. 
We write $\epsilon \colon A \to \Q$ for the
augmentation and   $A^{+}$ for the  augmentation ideal. When 
appropriate, we will
view $\Q$ as the DG algebra concentrated in degree 
$0$ with trivial differential and $\epsilon$ as a DG 
algebra map. 
A nilpotent, finite type CW complex 
$X$ admits a Sullivan minimal model
$\mathcal{M}_X, d_X$ which is a minimal DG algebra. 
Writing $\mathcal{M}_X = \Lambda V$ for some   graded vector space 
$V$, we recall that if $X$ is a {\em simple space} (that is,  the 
fundamental group of $X$ is abelian and acts trivially on the homotopy groups of 
$X$) then the rational homotopy groups of $X$ are recovered
by $V.$ Specifically, given graded spaces $V$ and $W$, let   $\Hom_n(V, W)$   
denote the space of linear maps between the graded spaces $V$ and $W$ 
reducing degrees by $n$. We then have Sullivan's isomorphism 
$\pi_n(X) \otimes \Q \cong \Hom_n(V, \Q)$
\cite[Th.15.11]{F-H-T}.  
A map of spaces $f \colon X \to Y$ has a minimal model which
is a map of DG algebras
$\mathcal{M}_f \colon \mathcal{M}_Y \to \mathcal{M}_X.$

Let  $A,d_A$ and $B,d_B$
be   DG algebras and
 $\phi \colon A \rightarrow B$  a DG algebra map.    
We say  $\theta \in \Hom_n(A, B)$
is a  \emph{$\phi$-derivation}
of degree $n$ if  $\theta (xy) = \theta(x)\phi(y) -
(-1)^{n|x|} \phi(x) \theta(y)$.  Let $\Der_{n}(A, B;
{\phi})$ denote the vector space of $\phi$-derivations of degree
$n$, for $n\geq 0$. Define a linear map $D_\phi \colon
\Der_{n}(A, B; {\phi}) \to \Der_{n-1}(A, B; {\phi})$ by
$D_\phi(\theta) = d_B \circ \theta  - (-1)^{|\theta|} \theta \circ
d_A$.   Then
$\Der_{*}(A, B; {\phi}), D_\phi$ is a chain complex. 
We write $H_n\big(\Der(A, B;
{\phi})\big)$ for the homology in degree $n$. 

Given a map $f \colon X \to Y,$   let $\map(X, Y; f)$ denote the 
path component of $f$ in the space of (unbased) continuous functions from $X$ to $Y.$ 
When $X$ and $Y$ are simply connected with $X$  a finite complex, 
 we have $ \pi_n(\map(X, Y;f)) \otimes \Q \cong 
 H_n(\Der(\mathcal{M}_Y, \mathcal{M}_X; \mathcal{M}_f)) $
for all $n \geq 2$  \cite[Th.2.1]{L-S1}.   We  recall  this
identification and 
describe a minor extension here. 

Given  $F \colon B \to \map(X, Y;f)$  
with  $B$ a simple CW complex of finite type,  
we describe the map induced by $F$ on rational homotopy groups. 
Let $\mathcal{F} \colon B \times X  \to Y$ denote  the adjoint map so 
that   $\mathcal{F} \circ i_2 
\sim f$ where $i_2 \colon X \to B \times X$ is the inclusion. 
Write    $\mathcal{M}_B = 
\Lambda W$ and $\mathcal{M}_X = \Lambda V.$  
Then  
  $\mathcal{M}_{\mathcal{F}}  \colon \mathcal{M}_Y \to \mathcal{M}_{B \times X} 
= \Lambda(W \oplus V)$.  Fix a homogeneous
basis $W = \Q(w_1, w_2, \ldots).$   
Given $\chi \in \mathcal{M}_Y$, we may then write
$ \mathcal{M}_{\mathcal{F}}(\chi) = \mathcal{M}_f(\chi) + \sum_{j} 
w_j \psi_j(\chi)  + A_\mathcal{F}(\chi)$
where $\psi_j \in \Hom_{|w_j|}(\mathcal{M}_Y, \mathcal{M}_X)$ 
and  $A_\mathcal{F}(\chi)$ is in the ideal of $\Lambda(W \oplus V)$ generated by 
the decomposables of $\Lambda W$.  
A standard check shows that each $\psi_j$    is an 
$\mathcal{M}_f$-derivation and a  cycle.   
We define 
$$ \Psi_F \colon \Hom_*(W, \Q) \to H_*(\Der(\mathcal{M}_Y, 
\mathcal{M}_X; \mathcal{M}_f))$$ by setting 
$\Psi_F(w_j^*) = \langle \psi_j \rangle$ and extending
by linearity. 
Here $w_j^* \in \Hom_{|w_j|}(W, \Q)$ denotes the dual of the basis element 
$w_j$.   
\begin{theorem}
\label{thm:adjoint}
Let $f \colon X \to Y$ be a map between simply connected CW complexes 
of finite type
with $X$ finite.  Let $F \colon B \to \map(X, Y;f)$ be a given map 
with $B$ a connected simple CW complex of finite type.  Then, with notation as above, we have
commutative diagrams
$$\xymatrix{    
\pi_n(B) \otimes \Q \ar[rr]^{\! \! \! \! \! (F_{\#})_\Q} 
\ar[d]_{\cong}^{\Phi_B} && 
\pi_n(\map(X, Y;f)) \otimes \Q \ar[d]^{\cong}_{\Phi_f}  \\
\Hom_n(W, \Q) \ar[rr]^{ \! \! \! \! \! \Psi_F} && H_n(\Der(\mathcal{M}_Y, 
\mathcal{M}_X; \mathcal{M}_f))}
$$
for all $n \geq 2$.  In the case $X = Y$ and $f = 1$ the result holds 
for $n = 1$ as well.
\end{theorem}
\begin{proof}
The map $\Phi_f$   is defined in the proof of \cite[Th.2.1]{L-S1}
 by following   the procedure above 
with $B = S^n$ but replacing  the minimal model of $S^n$ with its 
rational cohomology. The map $\Phi_B$ 
corresponds to the case $B = S^n, X = *$ and $Y= B$ is readily seen to 
be Sullivan's isomorphism, as recalled above.  
Compatibility with $\Psi_F$ is thus direct from 
definitions. For $n \geq 2$ the fact that $\Phi_f$ is an isomorphism 
is  \cite[Th.2.1]{L-S1}. 
When $n=1$, the map $\Phi_f$   extends  to a well-defined set 
map and is    a surjection \cite[Th.2.1c]{L-S2}, but     
will not generally be  a homomorphism  in degree $1$.
When $Y = X$ and $f = 1$, however,  
  the multiplication in $\pi_1(\map(X, X; 1))$ is induced by the 
  multiplication in $\map(X, X ;1)$, which is an $H$-space under 
  composition of maps. 
 Thus, the adjoint of the product of
two classes  
  $\alpha, \beta \colon S^1 \to \map(X, X;1)$ is given  
   by 
$$\xymatrix{ S^1   \times X  \ar[rr]^{ \Delta \times 1    } &&  S^1 \times 
S^1 \times X  \ar[rr]^{\Gamma} && X }
$$ 
where $\Gamma( z_1, z_2, x) = \left( \alpha(z_1) \circ \beta(z_2) 
\right) (x)$. Since $1 \colon X \to X$ is a two-sided identity in 
$\map(X, X;1),$ we have  
$\Gamma \circ (i_1, 1) = A$ and  $\Gamma \circ (i_2, 1) = B,$ where 
$A$ and $B$ denote the adjoints to $\alpha$ and $\beta,$ respectively, 
and $i_1, i_2 \colon S^1 \to S^1 \times S^1$ are the inclusions.  It 
follows easily from this that    $\Phi_f$ is a 
homomorphism in degree $1$.     Injectivity then follows by  the argument in \cite[Th.2.1]{L-S1} which only
requires $\Phi_f$ a homomorphism to extend to $n = 1$. 
\end{proof}

Now fix a fibration   $\xi \colon X \stackrel{j}{\to} E
\stackrel{p}{\to} B$  of simply connected finite type CW complexes   
with $X$ a finite
complex and  with classifying map    $h \colon B \to \Baut_1(X)$. 
We are interested in describing  $h$ at the level of rational homotopy 
groups  and so we may consider 
$\Omega h \colon \Omega B \to \Omega 
\Baut_1(X)$.  Using the equivalence
$\partial_U   \colon \Omega \Baut_1(X) \to \aut_1(X) = \map(X, X; 1)$ we obtain a map
\begin{equation} \label{eq:H}  H =   \partial_U \circ \Omega h \colon \Omega B \to 
\map(X, X; 1) \end{equation}
which fits the setting of \thmref{thm:adjoint}. 
Recall  that the Koszul-Sullivan model of the  fibration $\xi \colon X \stackrel{j}{\to} E
\stackrel{p}{\to} B$  is a short  exact sequence
\begin{equation}
\label{eq:KS}
\xymatrix{  \Lambda W, d_B \ar[r]^{\! \! \! \! \! \! \! \! \! P} &  
\Lambda(W \oplus V), d_E
\ar[r]^{\ \ \ \ \ \ J} & \Lambda V, d_X}\end{equation}
of  DG algebras.
The differential $d_E$ satisfies
$d_E(w) = d_B(w)$ for $w \in W$ and $d_E(v) -d_X(v) \in (\Lambda
W)^{+}\cdot \Lambda(W \oplus V)$ \cite[Prop.15.5]{F-H-T}.  The map $P$ is the inclusion
and the map   $J$ satisfies $J(w) =\epsilon(w)$ and
$J(v) = v.$
The DG algebras $\Lambda W, d_B$ and $\Lambda V, d_X$ are Sullivan
minimal models for $B$ and $X$ respectively.
The DG algebra $\Lambda(W \oplus V), d_E$
is a Sullivan model for the total space $E$ but is not, in general, a
minimal DG algebra.
Given $\chi \in \Lambda V$  we may then write
$d_E(\chi) = d_X (\chi) + \sum_j w_j \theta_{j}(\chi)  + B_E(\chi)$
where  $B_E(\chi)$ is in the ideal   of
$\Lambda(W \oplus V)$ generated by the decomposables in $\Lambda W$. 
Again we check directly that $\theta_j \in
\Der_{|w_j|-1}(\Lambda V, \Lambda V ; 1)$  and is a cycle. 
Define
$$\Theta_{\xi} : \Hom_n(W, \Q) \to H_{n-1}(\Der(\Lambda V, \Lambda 
V ; 1))$$
by setting
  $\Theta_{\xi}(w_j^*) = \langle \theta_j \rangle$  and extending
by linearity.  

\begin{theorem}
\label{thm:h} 
Let $\xi \colon X \stackrel{j}{\to} E
\stackrel{p}{\to} B$ be a fibration of simply connected CW complexes 
of finite type with $X$ finite and classifying map $h \colon B \to \Baut_1(X).$  
Then, with notation as above, we have a
commutative diagram
$$\xymatrix{  \pi_n(B) \otimes \Q \ar[rr]^{\! \! \! \! \! (h_{\#})_\Q } 
\ar[d]^{\Phi_B}_{\cong} && 
\pi_{n}(\Baut_1(X)) \otimes \Q \ar[d]_{\Phi_{1_X}'}^{\cong}  \\
\Hom_n(W, \Q) \ar[rr]^{ \! \! \! \! \! \Theta_\xi} && 
H_{n-1}(\Der(\Lambda V, \Lambda V; 1))}
$$ for $n \geq 2$. 
\end{theorem}
\begin{proof}
The map $\Phi'_{1_X}$   is obtained by precomposing $\Phi_{1_X}$ from 
\thmref{thm:adjoint}  with the identification $\pi_{n}(\Baut_1(X)) 
\otimes \Q \cong \pi_{n-1}(\map(X, X; 
1)) \otimes \Q.$ 
Let $\overline{W}$ denote the desuspension of $W$   
so that $\pi_*(\Omega B) \otimes \Q \cong \Hom_*(\overline{W}, \Q).$ Define 
$\overline{\Theta}_\xi \colon \Hom_*(\overline{W}, \Q) \to 
H_*(\Der(\Lambda V, \Lambda V; 1))$ by 
$\overline{\Theta}_\xi(\overline{w}_j^*) = \Theta_\xi(w_j^*).$  We 
show $\Psi_H = \overline{\Theta}_\xi$; the result follows. 

We first argue that the adjoint $\mathcal{H} \colon \Omega B \times X
\to X$  of the map  $H$ in  (\ref{eq:H}) is homotopic to the  {\em holonomy 
action} $\mathcal{H}_\xi = s_0 \circ i_0$  of 
 $\xi$ as defined by the following diagram:
 $$ \xymatrix{& \Omega B \times X  \ar[d]^{i_0} \ar[rd]^{i}
 \ar[d]^{i_0} \ar@/_4pc/[dd]_-{\mathcal{H}_\xi} \\
 &  PB \times_B E  \ar@/_.25pc/[d]_{s_0} \ar[r]^{i'} & MB \times_B E \ar[d]^{\pi_2}
  \ar[r]^{\ \ \pi_1} \ar[dr]^{q_1} & MB \ar[d]^{q_0} \\
 & X \ar@/_.25pc/[u]_{j_0}\ar[r]^j & E \ar[r]^p & B.}
 $$
 Here, with notation as in  \cite[\S 2]{F-H-T},  $MB$ is the
 space of Moore paths on $B$ with $q_0, q_1$ evaluation at $0$ and the
 length of the path, respectively,   $PB$   the
 subspace  of paths
 which end    at the basepoint of $B$ and 
 the other maps are the 
 the obvious inclusions and   projections. The inclusion $j_0$ is a homotopy
 equivalence \cite[Prop.2.5(ii)]{F-H-T} and we have denoted a homotopy
 inverse by $s_0$.  
Let $\mathcal{H}_\infty \colon  \Omega \Baut_1(X) \times  X    \to
  X$ denote the holonomy action of the universal
  $X$-fibration. Then $\mathcal{H}_\xi \sim \mathcal{H}_\infty \circ 
  \left(  \Omega h \times 1_X  \right)$  by the naturality 
of the holonomy action with respect to pull-backs
\cite[Prop.11.4]{Hilton}. Taking adjoints, we obtain $H_\xi \sim
  H_\infty \circ \Omega h.$   Finally,  
  $H_\infty \sim \partial_U \colon \Omega \Baut_1(X) \to \aut_1(X)$ by Stasheff's formulation of the 
  classification theorem  in terms of parallel transports \cite{St2}:
   the maps $\partial_U$  and $H_\infty$
  are transports giving the same principal $\aut_1(X)$-fibration
  (namely, the universal one) and hence are homotopic.

On \cite[p.419]{F-H-T}, the authors obtain a  Sullivan  model 
for the diagram $\Omega B \times X    \stackrel{i_0}{\lto}  
PB \times_{B} E
\stackrel{j_0}{\longleftarrow} X$ occuring in the definition of the
holonomy action $\mathcal{H}_\xi$.  In our notation, this is a diagram
of DG algebras
$$ \Lambda(\overline{W} \oplus V   ), \overline{d}_X \stackrel{I}{\longleftarrow}
\Lambda(W  \oplus \overline{W} \oplus V), \overline{d}_E
\stackrel{J}{\lto} \Lambda V,
d_X.$$
The differential $\overline{d}_X$ is given by $\overline{d}_X(v) = 
d_X(v)$   while
$\overline{d}_X(\overline{W}) = 0.$
The differential $\overline{d}_E$ satisfies  $\overline{d}_E(x) = d_E(x)$
for $x \in W \oplus V$.  The differential $\overline{d}_E$ on 
$\overline{w} \in \overline{W}$ is determined by the (push-out) 
construction involved.  Tracing this through, we see 
$\overline{d}_E(\overline{w}) - w \in \mathcal{J}$ where   
$\mathcal{J}$ is the ideal
of $\Lambda (W  \oplus \overline{W} \oplus V)$  generated by decomposables 
in $\Lambda(W \oplus \overline{W})$.  The maps $I$ and $J$ are projections.
A DG algebra map $S \colon \Lambda V, d_X \to \Lambda(W 
\oplus \overline{W} \oplus V), \overline{d}_E$ with $S \circ J = 1_{\Lambda V}$
exists by the standard lifting lemma for minimal models (cf. 
\cite[Lem.12.4]{F-H-T}). 
The map $S$ is a Sullivan model for the equivalence $s_0$ and  $I \circ 
S$ is one for  
$\mathcal{H}.$  Given $\chi  \in \Lambda V$ write
$S(\chi) = \chi + \sum_{j} \overline{w}_{j} \overline{\theta}_j(\chi) 
+  \sum_{j}w_j\varphi_j(\chi)    + C_S(\chi)$
for  $C_S(\chi) \in \mathcal{J}$. 
Using the fact that $S$ is a map of algebras, we obtain that
$\overline{\theta}_j \in \Der_{|w_j|-1}(\Lambda V, 
\Lambda V; 1) $ and $\varphi_j  \in \Der_{|w_j|}(\Lambda V, 
\Lambda V; 1) .$
Using the fact that $S$ is a chain map and
$\overline{d}_E(\mathcal{J}) \subseteq \mathcal{J},$ we obtain
$D_{1_X}(\overline{\theta}_j) = 0$ while
$D_{1_X}(\varphi_j) = \theta_j - 
\overline{\theta}_j.$
Thus $\Psi_{H}(\overline{w}_j^*) = \langle 
\overline{\theta}_j \rangle = \langle \theta_j 
\rangle = \overline{\Theta}_\xi(\overline{w}_j^*)$.
\end{proof}

\begin{remark} \label{rem:holonomy}
The {\em holonomy representation} of $\xi$  is  the action of  the
homotopy Lie algebra $\pi_*(\Omega B) \otimes \Q$ on $H_*(X; \Q)$ 
induced by the holonomy action $H_\xi \colon \Omega B \times X \to X$ 
(cf. \cite[p.415]{F-H-T}).
Define the  ``induced derivation'' map 
$I \colon H_*(\Der(\Lambda V, \Lambda V; 1)) \to \Der_*(H^*(X;\Q), 
H^*(X; \Q); 1)$ 
by setting
$I(\langle\theta\rangle)(\langle\chi\rangle) = \langle \theta(\chi)
\rangle$ for $\theta$ a derivation cycle and $\chi$ a cycle of
$\Lambda V, d_X.$ It is easy to check $I$ is well-defined.
By \thmref{thm:h} and \cite[Th.31.3]{F-H-T}, we see
$I \circ \Theta_\xi \colon \Hom_*(W, \Q) \to \Der_*(H^*(X;\Q), H^*(X; 
\Q) ; 1)$ is 
 dual to  the holonomy representation (up to sign). 
\end{remark}

\section{The Rationalized Gottlieb Sequence} \label{sec:rational}
The Gottlieb sequence  is   a $P$-local invariant
of the fibre-homotopy type of fibrations  of simply connected 
finite type CW 
complexes, provided the fibre is a finite complex. 
This fact    follows directly from     isomorphisms $G_n(X)
\otimes \Z_P \cong G_n(X_P)$ \cite{Lan} and $G_n(E, X; j) \otimes
\Z_P \cong G_n(E_P, X_P; j_P)$ \cite{Sm}.  
We define the $P$-local Gottlieb
homology $GH_n(\xi; \Z_P)$ of a fibration $\xi$ by setting
$GH_n(\xi; \Z_P) = GH_n(\xi_P).$ 
We will focus exclusively on the 
rational case here. We say a fibration $\xi$ is {\em rationally
Gottlieb trivial} if the sequence   (\ref{eq:Gottlieb sequence})
 splits into short exact sequences in each degree after tensoring with $\Q$.  

Consider, again, a fibration $\xi \colon X \stackrel{j}{\to} E 
\stackrel{p}{\to} B$   with $X$ a finite complex and classifying map 
$h \colon B \to \Baut_1(X)$. We wish to compare the condition that $\xi$ 
is rationally Gottlieb trivial,     with that of the vanishing of 
 $(h_\sharp)_\Q$. Since both of these conditions  
imply $\xi$  is {\em rationally weak-homotopically trivial}, 
that is, that   the   rationalized connecting homomorphism 
$(\partial_\sharp)_\Q = 0, $  we  may assume this 
condition holds for $\xi$. 
In this case, 
 the rationalization of the Gottlieb sequence for $\xi$ breaks up
into three-term sequences at each degree $n \geq 2$  of the form
\begin{equation} \xymatrix{
0 \ar[r] & G_n(X) \otimes \Q \ar[r]^{\! \! \! \! \! \! \! \! \! \! (j_\sharp)_\Q} &
G_n( E,  X; j) \otimes \Q  \ar[r]^{\ \ (p_\sharp)_\Q} & \pi_n(B) \otimes \Q
\ar[r] & 0.}
\label{eq:Got1}
\end{equation}
Moreover, the model of the total space of the fibration given in
(\ref{eq:KS}) for $\xi$  is a minimal DG algebra:
$\mathcal{M}_{E}, d_E = \Lambda(W \oplus V), d_E$ \cite[Prop.4.12]{Hal1}.
As above, we have  $d_E(w) =
d_B(w)$ while, for $\chi \in \Lambda V,$ $d_E(\chi) = d_X(\chi) + \sum_j w_j \theta_j(\chi) +
B_E(\chi)$ where $B_E(\chi)$ is  in the ideal $\mathcal{I}$ of $\Lambda(W \oplus V)$ generated 
by decomposables in $\Lambda W.$ 
The maps   $J \colon  \Lambda(W \oplus V ), d_E \to \Lambda V, d_X$ and
$P \colon  \Lambda W, d_B \to  \Lambda(W \oplus V ), d_E$
are Sullivan models for $j$ and $p$, respectively.  
They induce {\em linearization} maps $Q(J) \colon W \oplus V \to V $
and $Q(P) \colon W   \to W \oplus V$ (see \cite[p.171]{F-H-T}) which in 
turn induce  maps $Q(J)^* \colon \Hom_n(V, \Q)  \to \Hom_n(W \oplus V, \Q)$
and $Q(P)^* \colon \Hom_n(W \oplus V, \Q) \to \Hom_n(W, \Q)$ by 
taking duals. These maps correspond, under the identifications 
$\Phi_X, \Phi_E, \Phi_B$ of \thmref{thm:adjoint},  to the maps 
induced on rational homotopy groups by $j$ and $p.$

Composition with  $\epsilon \colon \Lambda V 
\to \Q$ 
induces a chain map $\epsilon_* \colon \Der_n(\Lambda(W \oplus V), 
\Lambda V; J) \to \Der_n(\Lambda (W \oplus V), \Q; \epsilon).$ 
The  minimality of $\mathcal{M}_E, d_E$ implies   $H_n(\Der(\Lambda (W \oplus V), \Q; 
\epsilon)) \cong \Hom_n(W \oplus V, \Q)$ and so we obtain a map
$$H(\epsilon_*) \colon
 H_n(\Der(\Lambda(W \oplus V), 
\Lambda V; J)) \to \Hom_n(\Lambda (W \oplus V), \Q).$$
Define the {\em $n$th rationalized evaluation subgroup  
of $J$} by $$G_n(\Lambda(W \oplus V), \Lambda V; J) = 
\im(H(\epsilon_*))$$ 
for $n \geq 2$. Thus  $w^* \in \Hom_n(W, \Q)$ is in the 
subspace $G_n(\Lambda(W \oplus V), \Lambda V; J)$ if and only if  $w^*$ extends to 
a $J$-derivation  $\theta$ of degree $n$ with $D_J(\theta) = 0.$ 
We write $G_n(\Lambda V) = G_n(\Lambda V, \Lambda 
V; 1)$ and call this the  {\em $n$th rationalized Gottlieb group} of 
$\Lambda V , d_X$.  We then obtain a sequence 
\begin{equation} \label{eq:Got}
 \xymatrix{ 0 \ar[r] & G_n(\Lambda V) \ar[r]^{\! \! \! \! \! \! \! \! \! \! \! \! 
 \! \! \! \! \! \! \! \! \! \! \! \! Q(J)^*} &
G_n(\Lambda(W \oplus V), \Lambda V; J) \ar[r]^{ \ \ \ \ \ \ \ \  Q(P)^*} & \Hom_n(W, \Q)
\ar[r] & 0 }
\end{equation}
for each $n \geq 2$. The following is a direct consequence of 
\cite[Cor.2.2]{L-S1}. 
\begin{theorem} \label{thm:Gseq} Let $\xi \colon X \stackrel{j}{\to} E \stackrel{p}{\to} B$ be a
fibration of simply connected, finite type CW complexes with $X$ 
finite and  $(\partial_\sharp)_\Q = 0.$ Then, for each $n \geq 2,$
the sequence {\em (\ref{eq:Got1})} is equivalent to the sequence  
{\em (\ref{eq:Got})}.
\qed
\end{theorem}

Using these identifications,  we extend
the result above for spherical fibrations (\thmref{thm:GH=0}) to the 
following:
\begin{theorem}
\label{thm:rational split}
Let $\xi : X \stackrel{j}{\rightarrow} E \stackrel{p}{\rightarrow} B$
be a fibration of simply connected, finite type  CW complexes   with
$X$ finite and classifying map $h \colon B \to \Baut_1(X)$. 
The following are equivalent:
\begin{itemize}
\item[(1)] $(h_\sharp)_\Q = 0 \colon \pi_n(B) \otimes \Q \to 
\pi_n(\Baut_1(X)) \otimes \Q$ for each $n \geq 2,$ 
\item[(2)]  $\xi$ is rationally Gottlieb  trivial. That is,    we have   short 
exact sequences
$$ \xymatrix{ 0 \ar[r] &
 G_n(X) \otimes \Q \ar[r]^{\! \! \! \! \! \! \!  (j_\sharp)_\Q} & G_n(E,
X; j) \otimes \Q \ar[r]^{\ \   (p_\sharp)_\Q} & \pi_n(B) \otimes \Q
\ar[r] &  0}$$ for each $n \geq 2,$
\item[(3)]  $(p_\sharp)_\Q \colon G_n(E, X; j) \otimes \Q \to 
\pi_n(B) \otimes \Q$ is surjective for all $n \geq 2.$
\end{itemize}

\end{theorem}
\begin{proof}
We begin with the implication $(1) \Rightarrow (2).$ 
By virtue of \thmref{thm:one-to-one}, we need only show exactness at
the middle term, that is, that $GH_*(\xi;\Q) = 0$. For this, suppose
that $\psi \in \Der_n(\Lambda(W\oplus V),\Lambda V; J)$ is a
$D_J$-cycle for which there exists $v \in V$ with $\psi(v) = 1$.
(This corresponds to an element $x \in G_n(E,X;j)\otimes\Q$ with
$(p_\sharp)_\Q(x) = 0$.)  We may also assume that $\psi(w_j) = 0$
for $|w_j| \leq n$.  Write $D$ for the differential in
$\Der_*(\Lambda V, \Lambda V;1)$.  We wish to find a $D$-cycle
$\theta \in \Der_n(\Lambda V, \Lambda V;1)$ that satisfies
$\theta(v) =1$. (This corresponds to the element $x \in
G_n(E,X;j)\otimes\Q$ being in the image of $G_n(X) \otimes \Q$ under $(j_\sharp)_\Q$.)  As a
first approximation to such a $\theta$, consider the derivation
$\psi_X \in \Der_n(\Lambda V, \Lambda V;1)$ obtained by simply
restricting $\psi$ to $\Lambda V$.  Since $D_J(\psi) = 0$, for $\chi
\in \Lambda V$ we have
$$\begin{aligned}
0 = D_J(\psi)(\chi) &= d_X \psi(\chi) - (-1)^{n}\psi d_E(\chi) \\
& = d_X \psi(\chi) - (-1)^{n}\psi d_X(\chi) - (-1)^{n}\psi\big(
\sum_j w_j \theta_j(\chi) + B_E(\chi) \big).
\end{aligned}
$$
Since $\psi$ is a $J$-derivation, and $J(W) = 0$, we have $\psi\big(
w_j \theta_j(\chi)\big) = \psi(w_j)\theta_j(\chi)$ and
$\psi\big(B_E(\chi)\big) = 0$.  This yields the identity
\begin{equation}\label{eq:psi_X obstruction}
0 = D(\psi_X)(\chi) - (-1)^n \sum_j \psi(w_j)\theta_j(\chi).
\end{equation}
The sum on the right-hand side is thus an obstruction to $\psi_X$
being a $D$-cycle, as we would wish.  With our hypothesis on the
classifying map, we may overcome this obstruction as follows.
Since $(h_\sharp)_\Q = 0$, \thmref{thm:h} implies that each
$\theta_j \in \Der_{|w_j|-1}(\Lambda V, \Lambda V;1)$ is a
$D$-boundary: $\theta_j = D(\varphi_j)$ for some
$\varphi_j \in \Der_{|w_j|}(\Lambda V, \Lambda V;1)$. For each
$j$, define derivations $\widehat{\theta_j} \in \Der_{n-1}(\Lambda V,
\Lambda V;1)$ and $\widehat{\varphi_j} \in \Der_{n}(\Lambda V,
\Lambda V;1)$ by setting
$$\widehat{\theta_j}(\chi) = \psi(w_j)\theta_j(\chi) \qquad
\text{and} \qquad \widehat{\varphi_j}(\chi) =
\psi(w_j)\varphi_j(\chi)$$
for $\chi \in \Lambda V$.  Since $\psi$ is a $J$-derivation and a
$D_J$-cycle, and since $J(W) = 0$, we see easily that $d_X\big(
\psi(w_j)\big) = 0$. Using this, a straightforward computation leads
to the identity
\begin{equation}\label{eq:hat derivations}
D(\widehat{\varphi_j}) = (-1)^{|w_j| - n} \widehat{\theta_j}.
\end{equation}
Now set $\theta = \psi_X - \sum_j (-1)^{|w_j|} \widehat{\varphi_j}
\in \Der_{n}(\Lambda V, \Lambda V;1)$. Since $\psi(w_j)$ is never
scalar, we have $\theta(v) = 1$. The identities (\ref{eq:psi_X
obstruction}) and (\ref{eq:hat derivations}) now give $D(\theta) =
0$, as required.

The implication $(2) \Rightarrow (3)$ is immediate. We prove $(3) 
\Rightarrow (1)$. 
By \thmref{thm:h}, it suffices to show that each of the derivations
$\theta_j \in \Der_{|w_j|}(\Lambda V, \Lambda V; 1)$ are
$D$-boundaries.  Let $\{w_j\}_{j \in J}$ be a
well-ordered homogeneous basis of $W$.  Because $(p_\sharp)_\Q$ is surjective,
each $w_j^* \in \Hom_{|w_j|}(W, \Q)$ extends to a $J$-derivation
$\eta_j \in \Der_{|w_j|}(\Lambda(W \oplus V), \Lambda V; J)$ with
$D_J(\eta_j) = 0$ and $\eta_j(w_j) = 1$. Without loss of generality,
we may suppose that $\eta_j(w_i) = 0$ for $i < j$ and that
$\eta_j(w_i)$ is of positive degree (or zero) for $i > j$.

For each $j$, let $\eta_{j,X} \in \Der_{|w_j|}(\Lambda V, \Lambda V;
1)$ denote the derivation obtained by restricting $\eta_j$ to
$\Lambda V$.  Then we have, for each $j$,
\begin{equation}\label{eq:D(psi_j)}
D\big( (-1)^{|w_j|} \eta_{j,X} \big) = \theta_j + \eta_j(w_{j+1})
\theta_{j+1} + \eta_j(w_{j+2})\theta_{j+2} + \cdots.
\end{equation}
This is simply a re-written version of (\ref{eq:psi_X obstruction}),
above.  Now choose a
particular $\theta_r$.  We will show that we may add a suitable
$D$-boundary to (\ref{eq:D(psi_j)}), so as to remove the terms from
the right-hand side other than $\theta_r$.  For this, we suppose
inductively that, for each $k \geq 1$, we have a derivation
$\widehat{\varphi_{k-1}} \in \Der_{|w_r|}(\Lambda V, \Lambda V; 1)$
such that
\begin{equation}\label{eq:4.3_ind_hyp}
D\big( (-1)^{|w_r|} \eta_{r,X} - \sum_{t = 1}^{k-1}
\widehat{\varphi_t} \big) = \theta_r + c_{r+k} \theta_{r+k} +
c_{r+k+1} \theta_{r+k+1} + \cdots,
\end{equation}
with each $c_{r+s} \in \Lambda V$, of degree $|c_{r+s}| = |w_{r+s}|
- |w_r|$, a sum of terms in $\Lambda ^{+} V$ each of the form
$$\eta_{r}(w_{p_1})\eta_{p_1}(w_{p_2}) \cdots
\eta_{p_l}(w_{r+s}),$$
for $r < p_1 < p_2 < \cdots < p_l < r+s$.

For the inductive step, define $\widehat{\varphi_{k}} = (-1)^{|w_r|}
c_{r+k} \eta_{r+k, X}$.  Since each $\eta_j$ is a $J$-derivation and a
$D_J$-cycle, and since $J(W) = 0$, we have that $d_X\big(
\eta_j(w_{i})\big) = 0$ for $i > j$.  It follows that $d_X(c_{r+s}) = 0$ for
each $s$.  Then
$$
\begin{aligned}
D\big( \widehat{\varphi_{k}} \big) &= (-1)^{|w_{r+k}|}c_{r+k} D( \eta_{r+k, X} )\\
&= c_{r+k} ( \theta_{r + k} + \eta_{r+k}(w_{r+k+1})\theta_{r+k+1} +
\eta_{r+k}(w_{r+k+2})\theta_{r+k+2} + \cdots), \end{aligned}
$$
from (\ref{eq:D(psi_j)}) with $j = r + k,$ and hence we have
$$
\begin{aligned}
D\big( (-1)^{|w_r|} \eta_{r,X} - \sum_{t = 1}^{k}
\widehat{\varphi_t} \big)
  = \theta_r & + \big( c_{r+k+1} - c_{r+k}
\eta_{r+k}(w_{r+k+1}) \big) \theta_{r+k+1}\\
&   + \big( c_{r+k+2} - c_{r+k} \eta_{r+k}(w_{r+k+2}) \big)
\theta_{r+k+2} + \cdots.
\end{aligned}
$$
The coefficients $c_{r+k+s} - c_{r+k} \eta_{r+k}(w_{r+k+s})$ of
each $\theta_{r+k+s}$ in this expression are of the required form,
and the inductive step is proven. Induction starts with $k = 1$,
where we have (\ref{eq:D(psi_j)}) for $j = r$ (take
$\widehat{\varphi_{0}} = 0$).

By induction, we may write
$$
D\big( (-1)^{|w_r|} \eta_{r,X} - \sum_{t \geq 1} \widehat{\varphi_t}
\big) = \theta_r.
$$
Notice that each $\widehat{\varphi_{t}}$ is a derivation, and that,
for any given $\chi \in \Lambda V$, we have
$\widehat{\varphi_{t}}(\chi) = 0$ for all but a finite number of
$t$.  Indeed, since $\eta_{r+k, X}$, used in the definition of
$\widehat{\varphi_{k}}$, decreases degree by $|w_{r+k}|$, we will
have $\widehat{\varphi_{t}}(\chi) = 0$ for all $t$ with $|w_{r+t}| \geq
|\chi|$.  The infinite sum is ``locally finite," therefore, and
defines a derivation.  This proves that each $\theta_r$ is a
$D$-boundary in $\Der_{|w_r|}(\Lambda V, \Lambda V; 1)$, as
required.
\end{proof}

 \thmref{thm:rational split}  provides a 
 link between the rationalized Gottlieb sequence and well-known 
 results on rational L-S 
 category.  Let $\cat_0(X)$ denote the rational L-S category of
 the space $X$, that is, $\cat_0(X) = \cat(X_\Q)$ where $X_\Q$ is the
 rationalization of $X$ and $\cat(X)$ denotes the ordinary L-S 
 category of $X$.
 \begin{corollary}
 \label{cor:cat}
 Let $\xi : X \stackrel{j}{\rightarrow} E \stackrel{p}{\rightarrow} B$
 be a fibration of simply connected, finite type  CW complexes   with
 $X$ finite. If $\xi$ satisfies one of the three equivalent 
 conditions of \thmref{thm:rational split},  then
 \begin{itemize}
 \item[(1)] 
 $\dim (G_{odd}(E, X; j) \otimes \Q) \leq \cat_0(X) +
  \dim(\pi_{odd}(B) \otimes \Q), $
 \item[(2)]
 $G_{even}(E, X; j) \otimes \Q \cong \pi_{even}(B) \otimes \Q.$
 \end{itemize}
If, further, $B$ has finite-dimensional rational homotopy then
\begin{itemize}
\item[(3)] $ \cat_0(E) \geq \cat_0(X) + \cat_0(B).$
\end{itemize}
 \end{corollary}
\begin{proof}
 The first two results are direct consequence of   
 \cite[Th.III]{F-H} and  condition (2) of \thmref{thm:rational split} .
 For the third result, observe that if  $(h_\sharp)_\Q = 0,$   
then by \remref{rem:holonomy} the   holonomy representation  of $\xi$ 
 is rationally trivial. The  result in this case follows from  
 \cite[Th.2]{Jess}.
\end{proof}

With rational techniques now at our disposal, we may rapidly add to
the store of illustrative examples begun in the integral setting
above.  The following provides a further example of a
weak-homotopically trivial fibration that is not Gottlieb trivial.

\begin{example} \label{ex:not Gottlieb}
Consider a fibration of the form $\xi \colon \C P^2
\stackrel{j}{\to} E  \stackrel{p}{\to} S^4$.  Any such fibration
admits a section and, in particular, is weak-homotopically trivial.
This follows by reasoning as in \exref{ex:integral not Gottlieb trivial}, using the fact
that $\pi_3(\C P^2) = 0$.  Now consider classifying maps for such a
$\xi$.  These are elements of $\pi_4\big(\Baut_1(\C P^2) \big) \cong
\pi_3\big(\map(\C P^2, \C P^2; 1)\big)$.  Let $\Lambda V, d_X =
\Lambda(v_2, v_5)$, with $d_{X}(v_2) =0$ and $d_{X}(v_5) = v_2^3,$
denote the minimal model of $\C P^2$.  Then the derivation $\theta
\in \Der_3(\Lambda V, \Lambda V; 1)$, defined by $\theta(v_2) = 0$
and $\theta(v_5) = v_2$, gives a non-zero class in
$H_3\big(\Der(\Lambda V, \Lambda V; 1) \big)$.  Hence,
$\pi_3\big(\map(\C P^2, \C P^2; 1)\big)\otimes \Q$ is non-zero and
so $\pi_4\big(\Baut_1(\C P^2) \big)$ contains elements of infinite
order.  Choose $h \colon S^4 \to \Baut_1(\C P^2)$ to be any map that
represents a class of infinite order. Since $(h_\#)_\Q \not= 0$, by  
\thmref{thm:rational split}  we have $(p_\#)_\Q \colon G_*(E,\C P^2;
j)\otimes \Q \to \pi_*(S^4)\otimes\Q$ is not surjective.  We conclude 
that
$p_\# \colon G_*(E,\C P^2; j) \to \pi_*(S^4)$ cannot be surjective.

It is perhaps interesting to note that the minimal model of $\xi$ is
determined by the requirement that its classifying map be
non-trivial rationally. Write $\mathcal{M}_{S^{4}} = \Lambda(w_4,
w_7)$ with $d_{B}(w_4) =0$ and $d_{B}(w_7) = w_4^2$.  Then any
fibration with base $S^4$ and fibre $\C P^2$ has minimal model of
the form $\Lambda(w_4, w_7), d_B \to \Lambda(w_4, w_7, v_2, v_5),
d_E \to \Lambda(v_2, v_5), d_X$. The only possible ``twisting" of
the differential $d_E$, for degree reasons, is of the form $d_E(v_5)
= v_2^3 + c w_4v_2$, for $c \in \Q$.  If $c = 0$, then the fibration
is rationally trivial; we have $(p_\#)_\Q \colon G_*(E,\C P^2;
j)\otimes \Q \to \pi_*(S^4)\otimes\Q$ surjective, as is easily
checked by the methods of this section; and hence by 
\thmref{thm:rational split} its classifying map is rationally trivial. If we assume this
is not the case, then the minimal model of the fibration is
determined, up to isomorphism, by the differential $d_E(v_5) = v_2^3
+ w_4v_2$.
\end{example}

%

Our final example illustrates   that the dimension of the  Gottlieb homology
and thus that of the  rationalized evaluation subgoup  of a fibre 
inclusion  can be
arbitrarily large   even when  the dimension of  $G_*(X)\otimes \Q $ and $\pi_*(B) \otimes
 \Q$  are small.

\begin{example}  \label{ex:GH=n} Let $n,k > 0$ be any    odd  integers. We
construct a rational fibration
$\xi_\Q : X_\Q \stackrel{j}{\rightarrow} E_\Q \stackrel{p}{\rightarrow} 
(S^{kn})_\Q$
where $X$ has the rational homotopy type of a finite complex with
  $\dim(G_*(X) \otimes \Q)  = 1$ such that
$\dim(GH_k(\xi_\Q)) = n.$
Let $V = \Q(v_1, \ldots, v_{n+1}, u)$ where each $v_j$ is of
degree $k$ and $|u| = k(n+1)-1.$ Define $d_X$ on $\Lambda V$ by setting
$d_X(v_i) = 0$ and $d_X(u) = v_1 v_2 \cdots v_{n+1}$.   Note that  $\Lambda
V, d_X$ is an elliptic   model  and so  may be
realized as a finite complex $X$.   Define a K.S. model
$$ \xymatrix{ \Lambda(w_{kn}), 0  \ar[r]^{ \! \! \!
\! \! \! \! \! \! \! \!
\! \! \! \! \! \! \! \! P} &  \Lambda(  \Q(w_{kn}) \oplus
V), d_E \ar[r]^{  \ \ \ \ \ \ \
J}&  \Lambda V, d_X}$$
by setting
$d_E(u) = v_1v_2\cdots v_{n+1} + w_{kn} v_{n+1}$  with $d_E(w_{kn}) 
= d_E(v_i) = 0$ for $i =1, \ldots, n.$
We show each $v_i^* \in \Hom_k(V, \Q)$ extends to a $J$-derivation
cycle for $i= 1, \ldots, n.$  Define $\psi_i \in \Der(\Lambda(  \Q(w_{kn}) \oplus
V), \Lambda V;J)$ by setting $\psi_i(v_i) = 1,  \psi_i(v_j) = \psi_i(u) = 0$ for $j
\neq i$ and $\psi_i(w_{kn}) = -v_1v_2 \cdots \widehat{v_i} \cdots v_n.$
It is easy to check $D_J(\psi_i) = 0$ as needed.  Visibly $G_j(X)
\otimes \Q  = 0$  for $j \neq k(n+1) -1$ while $G_{k(n+1) -1}(X) 
\otimes \Q \cong \Q$.
Thus $\dim(GH_k(\xi_\Q)) = n$.

\end{example}
%
\providecommand{\bysame}{\leavevmode\hbox to3em{\hrulefill}\thinspace}
\providecommand{\MR}{\relax\ifhmode\unskip\space\fi MR }
\providecommand{\MRhref}[2]{%
  \href{http://www.ams.org/mathscinet-getitem?mr=#1}{#2}
}
\providecommand{\href}[2]{#2}

\end{document}